\documentclass[11pt]{article}
\usepackage{epic,latexsym,amssymb}
\usepackage{color}
\usepackage{tikz}

\newtheorem{thm}{Theorem}

\newtheorem{cor}{Corollary}

\newcommand{\qed}{$\Box$}

\newcommand{\smallqed}{{\tiny ($\Box$)}}

\newcommand{\D}{{Dominator }}
\newcommand{\St}{{Staller }}
\newcommand{\w}{{\rm w}}

\newcommand{\gtd}{\gamma_{{\rm tg}}}

\newenvironment{unnumbered}[1]{\trivlist \item [\hskip \labelsep {\bf
#1}]\ignorespaces\it}{\endtrivlist}

\newcommand{\proof}{\noindent\textbf{Proof. }}

\begin{document}

\title{Progress Towards the Total Domination Game $\frac{3}{4}$-Conjecture}

\author{$^1$Michael A. Henning and $^2$Douglas F. Rall
\\ \\
$^1$Department of Pure and Applied Mathematics \\
University of Johannesburg \\
Auckland Park, 2006 South Africa\\
\small \tt Email: mahenning@uj.ac.za  \\
\\
$^2$Department of Mathematics \\
Furman University \\
Greenville, SC, USA\\
\small \tt Email: doug.rall@furman.edu}

\date{}
\maketitle

\begin{abstract}
In this paper, we continue the study of the total domination game in graphs introduced in [Graphs Combin. 31(5) (2015), 1453--1462], where the players Dominator and Staller alternately select vertices of $G$. Each vertex chosen must strictly increase the number of vertices totally dominated, where a vertex totally dominates another vertex if they are neighbors. This process eventually produces a total dominating set $S$ of $G$ in which every vertex is totally dominated by a vertex in $S$. Dominator wishes to minimize the number of vertices chosen, while Staller wishes to maximize it. The game total domination number, $\gamma_{\rm tg}(G)$, of $G$ is the number of vertices chosen when Dominator starts the game and both players play optimally. Henning,  Klav\v{z}ar and Rall [Combinatorica, to appear] posted the $\frac{3}{4}$-Game Total Domination  Conjecture that states that if $G$ is a graph on $n$ vertices in which every component contains at least three vertices, then $\gamma_{\rm tg}(G) \le \frac{3}{4}n$. In this paper, we prove this conjecture over the class of graphs $G$ that satisfy both the condition that the degree sum of adjacent vertices in $G$ is at least~$4$ and the condition that no two vertices of degree~$1$ are at distance~$4$ apart in $G$. In particular, we prove that by adopting a greedy strategy, Dominator can complete the total domination game played in a graph with minimum degree at least~$2$ in at most~$3n/4$ moves.
\end{abstract}

{\small \textbf{Keywords:} Total domination game; Game total domination number; $3/4$-Conjecture } \\
\indent {\small \textbf{AMS subject classification:} 05C65, 05C69}

\newpage
\section{Introduction}

The domination game in graphs was first introduced by Bre\v{s}ar, Klav\v{z}ar, and Rall~\cite{BKR10} and extensively studied afterwards
in~\cite{BDKK14, brklra-2013, BKKR13, BKR13, Buj14, BuTu15, bujtas-2015+, doko-2013, HeLo-2015+, KWZ13, Kos} and elsewhere. A vertex \emph{dominates} itself and its neighbors. A \emph{dominating set} of $G$ is a set $S$ of vertices of $G$ such that every vertex in $G$ is dominated by a vertex in $S$.  The \emph{domination game} played on a graph $G$ consists of two players, \emph{Dominator} and \emph{Staller}, who take turns choosing a vertex from $G$. Each vertex chosen must dominate at least one vertex not dominated by the vertices previously chosen. The game ends when the set of vertices chosen becomes a dominating set in $G$. Dominator wishes to minimize the number of vertices chosen, while Staller wishes to  end the game with as many vertices chosen as possible. The \emph{game domination number},  $\gamma_g(G)$, of $G$ is the number of vertices chosen when Dominator starts the game and both players play optimally.

Much interest in the domination game arose from the $3/5$-Game Domination Conjecture posted by Kinnersley, West, and Zamani in~\cite{KWZ13}, which states that if $G$ is an isolate-free forest on $n$ vertices, then $\gamma_{g}(G) \le \frac{3}{5}n$. This conjecture remains open, although to date it is shown to be true for graphs with minimum degree at least~$2$ (see,~\cite{HeKi15}), and for isolate-free forests in which no two leaves are at distance~$4$ apart (see,~\cite{Buj14}).

Recently, the total version of the domination game was investigated in~\cite{hkr-2015}, where it was demonstrated that these two versions differ significantly. A vertex \emph{totally dominates} another vertex if they are neighbors. A \emph{total dominating set} of a graph $G$ is a set $S$ of vertices such that every vertex of $G$ is totally dominated by a vertex in
$S$. The \emph{total domination game} consists of two players called \emph{Dominator} and \emph{Staller}, who take turns choosing a vertex from~$G$. Each vertex chosen must totally dominate at least one vertex not totally dominated by the set of vertices previously chosen. Following the notation of~\cite{hkr-2015}, we call such a chosen vertex a \emph{legal move} or a \emph{playable vertex} in the total domination game. The game ends when the set of vertices chosen is a total dominating set in $G$. Dominator's objective is to minimize the number of vertices chosen, while Staller's is to end the game with as many vertices chosen as possible.

The \emph{game total domination number}, $\gtd(G)$, of $G$ is the number of vertices chosen when Dominator starts the game and both players employ a strategy that achieves their objective. If \St starts the game, the resulting number of vertices chosen is the \emph{Staller-start game total domination number}, $\gtd'(G)$, of $G$.

A \emph{partially total dominated graph} is a graph together with a declaration that some vertices are already totally dominated; that is, they need not be totally dominated in the rest of the game. In~\cite{hkr-2015}, the authors present a key lemma, named the \emph{Total Continuation Principle}, which in particular implies that when the game is played on a partially total dominated graph $G$, the numbers $\gtd(G)$ and $\gtd'(G)$ can differ by at most~$1$.

Determining the exact value of $\gtd(G)$ and $\gtd'(G)$ is a challenging problem, and is currently known only for paths and cycles~\cite{dh-2015+}. Much attention has therefore focused on obtaining upper bounds on the game total domination number in terms of the order of the graph. The best general upper bound to date on the game total domination number for general graphs is established in~\cite{hkr-2015+}.

\begin{thm}{\rm (\cite{hkr-2015+})}
If $G$ is a graph on $n$ vertices in which every component contains at least three vertices, then $\gtd(G) \le \frac{4}{5}n$.
\label{t:known}
\end{thm}

Our focus in the present paper is the following conjecture posted by Henning,  Klav\v{z}ar and Rall~\cite{hkr-2015+}.

\begin{unnumbered}{$\frac{3}{4}$-Game Total Domination Conjecture (\cite{hkr-2015+})}
\label{c:34}
If $G$ is a graph on $n$ vertices in which every component contains at least three vertices, then $\gtd(G) \le \frac{3}{4}n$.
\end{unnumbered}

%
Bujt\'{a}s, Henning, and Tuza~\cite{BuHeTu15} recently  proved the $\frac{3}{4}$-Conjecture over the class of graphs with minimum degree at least~$2$. To do this, they raise the problem to a higher level by introducing a transversal game in hypergraphs, and establish a tight upper bound on the game transversal number of a hypergraph with all edges of size at least~$2$ in terms of its order and size. As an application of this result, they prove that if $G$ is a graph on $n$ vertices with minimum degree at least~$2$, then $\gtd(G) \le \frac{8}{11}n$, which validates the 
$\frac{3}{4}$-Game Total Domination Conjecture on graphs with minimum degree at least~$2$.

For notation and graph theory terminology not defined herein, we in general follow~\cite{HeYe_book}. We denote
the \emph{degree} of a vertex $v$ in a graph $G$ by $d_G(v)$, or simply by $d(v)$ if the graph $G$ is clear
from the context.
The minimum degree among the vertices of $G$ is denoted by $\delta(G)$.
A vertex of degree~$1$ is called a \emph{leaf} and its neighbor a \emph{support vertex}.
If $X$ and $Y$ are subsets of vertices in a graph $G$, then the set $X$ \emph{totally dominates} the set
$Y$ in $G$ if every vertex of $Y$ is adjacent to at least one vertex of $X$. In particular, if $X$ totally
dominates the vertex set of $G$, then $X$ is a total dominating set in $G$. For more information on
total domination in graphs see the recent book~\cite{HeYe_book}. Since an isolated vertex in a graph cannot be totally dominated by definition, all graphs considered will be without isolated vertices.  We also use the standard notation $[k] = \{1,\ldots,k\}$.

\section{Main Result}

In this paper we prove the following result. Its proof is given in Section~\ref{S:proof}.

\begin{thm}
The $\frac{3}{4}$-Game Total Domination Conjecture is true over the class of graphs $G$ that satisfy both conditions {\rm (a)} and {\rm (b)} below: \\
\indent {\rm (a)} The degree sum of adjacent vertices in $G$ is at least~$4$. \\
\indent {\rm (b)} No two leaves are at distance~$4$ apart in $G$.
\label{t:main}
\end{thm}

As a special case of Theorem~\ref{t:main}, the $\frac{3}{4}$-Game Total Domination Conjecture is valid on graphs with minimum degree at least~$2$.

\begin{cor}{\rm (\cite{BuHeTu15})}
The $\frac{3}{4}$-Game Total Domination Conjecture is true over the class of graphs with minimum degree at least~$2$.
\label{c:Tdom2}
\end{cor}

\section{Proof of Main Result}
\label{S:proof}

In this section, we give a proof of our main theorem, namely Theorem~\ref{t:main}. For this purpose, we adopt the approach of the authors in~\cite{hkr-2015+} and color the vertices of a graph with four colors that reflect four different types of vertices. More precisely, at any stage of the game, if $D$ denotes the set of vertices played to date where initially $D = \emptyset$, we define as in~\cite{hkr-2015+} a \emph{colored}-\emph{graph} with respect to the played vertices in the set $D$ as a graph in which every vertex is colored with one of four colors, namely white, green, blue, or red, according to the following rules. \\ [-1.75em]
\begin{itemize}
\item A vertex is colored \emph{white} if it is not totally dominated by $D$ and does not belong to $D$. \\ [-2em]
\item A vertex is colored \emph{green} if it is not totally dominated by $D$ but belongs to $D$. \\ [-2em]
\item A vertex is colored \emph{blue} if it is totally dominated by $D$ but has a neighbor not totally dominated by $D$. \\ [-2em]
\item A vertex is colored \emph{red} if it and all its neighbors are totally dominated by $D$.
\end{itemize}

As remarked in~\cite{hkr-2015+}, in a partially total dominated graph the only playable vertices are those that have a white or green neighbor since a played vertex must totally dominate at least one new vertex. In particular, no red or green vertex is playable. Further, as observed in~\cite{hkr-2015+}, once a vertex is colored red it
plays no role in the remainder of the game, and edges joining two blue vertices play no role in the game. Therefore, we may assume a partially total dominated graph contains no red vertices and has no edge joining two blue vertices. The resulting graph is called a \emph{residual graph}. We note that the degree of a white or green vertex in the residual graph remains unchanged from its degree in the original graph.

Where our approach in the current paper differs from that in~\cite{hkr-2015+} is twofold. First, we define two new colors in the colored-graph that may possibly be introduced as the game is played. Second, our assignment of weights to vertices of each color differs from the assignment in~\cite{hkr-2015+}. Here, we associate a weight with every vertex in the residual graph as follows:

\begin{center}
\begin{tabular}{|c|c|}  \hline
\emph{Color of vertex} & \emph{Weight of vertex} \\ \hline
white & $3$ \\
green & $2$  \\
blue & $1$  \\
red & $0$  \\ \hline
\end{tabular}
\end{center}
\begin{center}
\textbf{Table~1.} The weights of vertices according to their color.
\end{center}

We denote the \emph{weight} of a vertex $v$ in the residual graph $G$ by $\w(v)$. For a subset $S \subseteq V(G)$ of vertices of $G$, the \emph{weight} of $S$ is the sum of the weights of the vertices in $S$, denoted $\w(S)$. The \emph{weight} of $G$, denoted $\w(G)$, is the sum of the weights of the vertices in $G$; that is, $\w(G) = \w(V(G))$. We define the \emph{value} of a playable vertex as the decrease in weight resulting from playing that vertex.

We say that \D can \emph{achieve his $4$-target} if he can play a sequence of moves guaranteeing that on average the weight decrease resulting from each played vertex in the game is at least~$4$. In order to achieve his $4$-target, \D must guarantee that a sequence of moves $m_1, \ldots, m_k$ are played, starting with his first move $m_1$, and with moves alternating between \D and \St such that if $\w_i$ denotes the decrease in weight after move $m_i$ is played, then
\begin{equation}
\sum_{i = 1}^k \w_i \ge 4k\,, \label{Eq1}
\end{equation}
and the game is completed after move $m_k$.
In the discussion that follows, we analyse how \D can \emph{achieve his $4$-target}.
For this purpose, we describe a move that we call a greedy move.
\begin{itemize}
\item A \emph{greedy move} is a move that decreases the weight by as much as possible. We say that \D follows a \emph{greedy strategy} if he plays a greedy move on each turn.
 \end{itemize}

We are now in a position to prove our main result, namely Theorem~\ref{t:main}. Recall its statement.

\noindent \textbf{Theorem~\ref{t:main}}. \emph{The $\frac{3}{4}$-Game Total Domination Conjecture is true over the class of graphs $G$ that satisfy both conditions {\rm (a)} and {\rm (b)} below: \\
\indent {\rm (a)} The degree sum of adjacent vertices in $G$ is at least~$4$. \\
\indent {\rm (b)} No two leaves are at distance~$4$ apart in $G$.
}

\medskip
\proof  Let $G$ be a graph that satisfies both conditions (a) and (b) in the statement of the theorem. Coloring the vertices of $G$ with the color white we produce a colored-graph in which every vertex is colored white. In particular, we note that $G$ has $n$ white vertices and has weight~$\w(G) = 3n$. Before any move of Dominator, the game is in one of the following two phases. \\ [-1.75em]
\begin{itemize}
\item \emph{Phase}~1, if there exists a legal move of value at least~$5$. \\ [-2em]
\item \emph{Phase}~2, if every legal move has value at most~$4$.
\end{itemize}
We proceed with the following claims.

\begin{unnumbered}{Claim~\ref{t:main}.1}
Every legal move in a residual graph decreases the total weight by at least~$3$.
\end{unnumbered}
\proof Every legal move in a colored-graph is a white vertex with at least one white neighbor or a blue vertex with at least one white or green neighbor. Let $v$ be a legal move in a residual graph. Suppose that $v$ is a white vertex, and so $v$ has at least one white neighbor. When $v$ is played, the vertex $v$ is recolored green while each white neighbor of $v$ is recolored blue, implying that the weight decrease resulting from playing~$v$ is at least~$3$. Suppose that $v$ is a blue vertex, and so each neighbor of $v$ is colored white or green. Playing the vertex $v$ recolors each white neighbor of $v$ blue or red and recolors each green neighbor of $v$ red. The weight of each neighbor of the blue vertex $v$ is therefore decreased by at least~$2$ when $v$ is played, while the vertex $v$ itself is recolored red and its weight decreases by~$1$. Hence, the total weight decrease resulting from playing~$v$ is at least~$3$.~\smallqed

\begin{unnumbered}{Claim~\ref{t:main}.2}
Let $R$ be the residual graph. If the game is in Phase~2 and if $C$ is an arbitrary component of $R$, then one of the following holds. \\
\indent {\rm (a)} $C \cong P_4$, with both leaves colored blue and both internal vertices colored white. \\
\indent {\rm (b)} $C \cong P_3$, with both leaves colored blue and with the central vertex colored green. \\
\indent {\rm (c)} $C \cong P_2$, with one leaf colored blue and the other colored green. \\
\indent {\rm (d)} $C \cong P_2$, with one leaf colored blue and the other colored white.
\end{unnumbered}
\proof Suppose the game is in Phase~2. We show first that every white vertex has at most one white neighbor in the residual graph $R$. Suppose, to the contrary, that a white vertex $v$ has at least two white neighbors. When $v$ is played the weight decreases by at least~$1 + 2 \cdot 2 = 5$, since the vertex $v$ is recolored green while each white neighbor of $v$ is recolored blue. This contradicts the fact that every legal move decreases the weight by at most~$4$.

We show next that every blue vertex has degree~$1$ in the residual graph $R$. Suppose, to the contrary, that a blue vertex $v$ has degree at least~$2$ in $R$. Playing the vertex $v$ recolors each white neighbor of $v$ blue or red and recolors each green neighbor of $v$ red. Thus, playing the vertex $v$ decreases the weight of each of its neighbors by at least~$2$. In addition, the vertex $v$ is recolored red, and so its weight decreases by~$1$. Hence, the weight decrease resulting from playing~$v$ is at least~$5$, a contradiction.

Suppose that $R$ contains a green vertex, $v$. Each neighbor of $v$ is colored blue, and, by our earlier observations, is therefore a blue leaf.
If $v$ is a leaf, then the component containing $v$ is a path isomorphic to $P_2$ with one leaf colored blue and the other colored green, and therefore satisfies condition~(c) in the statement of the claim. Hence, we may assume that the (green) vertex $v$ has at least two neighbors in $R$. If $v$ has at least three neighbors in $R$, then since every neighbor of $v$ is a blue leaf, the weight decrease resulting from playing an arbitrary neighbor of $v$ is at least~$5$, noting that such a move recolors $v$ and all its neighbors red. This produces a contradiction. Therefore, $v$ has exactly two neighbors in $R$, implying that the component containing $v$ is a path isomorphic to $P_3$ with both leaves colored blue and with the central vertex, namely $v$, colored green, and therefore satisfies condition~(b) in the statement of the claim. Hence, we may assume that there is no green vertex, for otherwise the desired result holds.

Suppose that there is a white vertex, $u$, in the residual graph $R$. Suppose that $u$ has no white neighbor. By our earlier observations, every neighbor of $u$ is a blue leaf. Playing a neighbor of $u$ therefore recolors all the neighbors of $u$ from blue to red. Since the degree of a white vertex in the residual graph remains unchanged from its degree in the original graph, we note in particular that $d_G(u) = d_R(u)$. If $u$ is not a leaf in $G$, then playing a neighbor of $u$ decreases the weight by at least~$3 + d_R(u) \ge 5$, a contradiction. Hence, $u$ is a leaf, and the component containing $u$ is a path isomorphic to $P_2$ with one leaf colored blue and the other colored white. We may therefore assume that the vertex $u$ has exactly one white neighbor, for otherwise the component containing $u$ satisfies condition~(d) in the statement of the claim.

Let $x$ be the white neighbor of $u$. Every neighbor of $u$ different from $x$ is a blue leaf, and every neighbor of $x$ different from $u$ is a blue leaf. Suppose that $u$ or $x$, say $u$, is a leaf. Since the degree sum of adjacent vertices in $G$ is at least~$4$, and $d_G(u) = d_R(u)$, the vertex $x$ has degree at least~$3$. Playing the vertex $u$ recolors $u$ from white to green, recolors $x$ from white to blue, and recolors all neighbors of $x$ different from $u$ from blue to red. Hence, playing $u$ decreases the total weight by at least~$5$, a contradiction. Therefore, neither $u$ nor $x$ is a leaf, implying that both $u$ and $x$ have at least one blue leaf neighbor. Suppose that $u$ or $x$, say $u$, has degree at least~$3$. Playing the vertex $x$ recolors $x$ from white to green, recolors $u$ from white to blue, and recolors each neighbor of $u$ different from $x$ from blue to red, implying that the total weight decrease resulting from playing~$x$ is at least~$5$, a contradiction. Therefore, both $u$ and $x$ have degree~$2$. Thus, the component containing $u$ and $x$ is a path isomorphic to $P_4$ with both leaves colored blue and both internal vertices colored white, and therefore satisfies condition~(a) in the statement of the claim.  This completes the proof of Claim~\ref{t:main}.2.~\smallqed

\medskip
By Claim~\ref{t:main}.2, once the game enters Phase~2 the residual graph is determined and each component satisfies one of the conditions~(a)--(d) in the statement of the claim.

\begin{unnumbered}{Claim~\ref{t:main}.3}
If the minimum degree in $G$ is at least~$2$, then \D can achieve his $4$-target by following a greedy strategy.
\end{unnumbered}
\proof Suppose that $\delta(G) \ge 2$ and \D follows a greedy strategy. Thus, at each stage of the game, \D plays a (greedy) move that decreases the weight by as much as possible. By Claim~\ref{t:main}.1, every move of Staller's decreases the weight by at least~$3$. Hence, whenever \D plays a vertex that decreases the weight by at least~$5$, his move, together with Staller's response, decreases the weight by at least~$8$. Therefore, we may assume that at some stage the game enters Phase~2, for otherwise  Inequality~(\ref{Eq1}) is satisfied upon completion of the game and \D can achieve his $4$-target.

Suppose that the first $\ell$ moves of \D each decrease the weight by at least~$5$, and that his $(\ell + 1)$st move decreases the weight by at most~$4$. Thus, $\w(m_{2i - 1}) + \w(m_{2i}) \ge 8$ for $i \in [\ell]$, and $\w(m_{2\ell + 1}) \le 4$. Let $R$ denote the residual graph immediately after Staller plays her $\ell$th move, namely the move $m_{2\ell}$. Thus,
\[
\sum_{i = 1}^{2\ell} \w_i = \sum_{i = 1}^{\ell} (\w(m_{2i - 1}) + \w(m_{2i})) \ge 8 \cdot \ell = 4 \cdot (2\ell).
\]

Since $\delta(G) \ge 2$, we note that $R$ contains no green or white leaf. Hence, by Claim~\ref{t:main}.2, every component $C$ of $R$ satisfies $C \cong P_4$, with both leaves colored blue and both internal vertices colored white, or $C \cong P_3$, with both leaves colored blue and with the central vertex colored green. If $C \cong P_4$, then $\w(V(C)) = 8$ and exactly two additional moves are required to totally dominate the vertices $V(C)$, while if $C \cong P_3$, then $\w(V(C)) = 4$ and exactly one move is played in $C$ to totally dominate the vertices $V(C)$. Suppose that $R$ has $t$ components isomorphic to $P_4$ and $s$ components isomorphic to $P_3$. Thus, $2t + s$ additional moves are needed to complete the game once it enters Phase~2. Further, these remaining $2t + s$ moves satisfy
\[
\sum_{i = 2\ell + 1}^{2\ell + 2t + s} \w_i = 4 \cdot (2t + s).
\]
Hence,
\[
\sum_{i = 1}^{2\ell + 2t + s} \w_i  = \sum_{i = 1}^{2\ell} \w_i + \sum_{i = 2\ell + 1}^{2\ell + 2t + s} \w_i \ge 4 \cdot (2\ell + 2t + s),
\]

\noindent
and so Inequality~(\ref{Eq1}) is satisfied upon completion of the game. Thus, \D can achieve his $4$-target simply by following a greedy strategy. This completes the proof of Claim~\ref{t:main}.3.~\smallqed

\medskip
We now return to the proof of Theorem~\ref{t:main}. By Claim~\ref{t:main}.3, we may assume that $G$ contains at least one leaf, for otherwise \D can achieve his $4$-target (and he can do so by following a greedy strategy).


As the game is played, we introduce a new color, namely \emph{purple}, which we use to recolor certain white support vertices. A purple vertex will have the same properties of a white vertex, except that the weight of a purple vertex is~$4$. The idea behind the re-coloring is that the additional weight of~$1$ assigned to a purple vertex will represent a \emph{surplus weight} that we can ``bank" and withdraw later. To formally define the recoloring procedure, we introduce additional terminology.

Consider a residual graph $R$ that arises during the course of the game. Suppose that $uvwx$ is an induced path in $R$, where $v$, $w$ and $x$ are all white vertices, and where $w$ is a support vertex and $x$ a leaf in $R$. We note that the vertex $u$ is colored white or blue. Such a vertex $u$ turns out to be problematic for Dominator, and we call such a vertex a \emph{problematic vertex}. Further, we call the path $uvwx$ a \emph{problematic path} associated with~$u$, and we call $w$ a \emph{support vertex associated with~$u$}. 

Suppose that \St plays a problematic vertex, $u$. Suppose that there are exactly $k$ support vertices, say $w_1,\ldots,w_k$, associated with~$u$. We note that $k \ge 1$. For $i \in [k]$, let $uv_iw_ix_i$ be a problematic path associated with~$u$ that contains~$w_i$. Thus, $v_i$, $w_i$ and $x_i$ are all white vertices, $w_i$ is a support vertex, and $x_i$ a leaf in $R$. Since no two leaves are at distance~$4$ apart in $G$, we note that if $k \ge 2$, then $v_i \ne v_j$ for $1 \le i,j \le k$ and $i \ne j$.

Suppose first that $k \ge 2$. In this case, playing the problematic vertex, $u$, decreases the total weight by at least~$2k+1$, since $u$ is recolored from blue to red or from white to green, while each neighbor $v_i$, $i \in [k]$, of $u$ is recolored from white to blue. Thus, the current value of $u$ is at least~$2k+1$. We now discharge the value of $u$ as follows. We discharge a weight of~$k$ from the value of $u$ and add a weight of~$1$ to every support vertex $w_i$, $i \in [k]$. Thus, by playing $u$ the resulting decrease in total weight is the value of $u$ in $R$ minus~$k$, which is at least~$k + 1 \ge 3$. Further, the weight of each (white) support vertex $w_i$, $i \in [k]$, increases from~$3$ to~$4$. We now re-color each support vertex $w_i$, $i \in [k]$, from white to purple.

Suppose secondly that $k = 1$ and the value of $u$ is at least~$4$. In this case, we proceed exactly as before: We discharge a weight of~$k = 1$ from the value of $u$, add a weight of~$1$ to the support vertex $w_1$, and re-color $w_1$ from white to purple. Thus, by playing $u$ the resulting decrease in total weight is at least~$3$.

In both cases, we note that the new weight of $w_i$, $i \in [k]$, is~$4$. Thus, every newly created purple vertex is a support vertex and has weight~$4$. We define a purple vertex to have the identical properties of a white vertex, except that its weight is~$4$. Thus, a purple vertex is not totally dominated by the vertices played to date and has not yet been played.

We note that if \St plays a problematic vertex, $u$, whose current value is at least~$4$, then the above discharging argument recolors every support vertex associated with~$u$ from white to purple. Further, by playing $u$ the resulting decrease in total weight is at least~$3$, and the weight of each newly created purple vertex is~$4$. We state this formally as follows.

\begin{unnumbered}{Claim~\ref{t:main}.4}
If \St plays a problematic vertex whose current value is at least~$4$, then the resulting decrease in total weight is at least~$3$.
\end{unnumbered}

We note further that if \St plays a problematic vertex, $u$, whose current value is exactly~$3$, the two internal vertices of the problematic path associated with~$u$ are unique. In particular, the support vertex associated with $u$ is unique.

We now introduce an additional new color, namely \emph{indigo}, which we use to recolor certain purple vertices. An indigo vertex will have the same properties of a blue vertex, except that the weight of an indigo vertex is~$2$ (while the weight of a blue vertex is~$1$). The idea behind the re-coloring is that the additional weight of~$1$ assigned to an indigo vertex will represent a \emph{surplus weight} that, as before, we can ``bank" and withdraw later.

More formally, suppose that a white leaf, say $z$,  adjacent to a purple vertex, say $x$, is played. When the leaf $z$ is played, it changes color from white to green and its support neighbor, $x$, changes color from purple to blue (noting that a purple vertex has the same properties as a white vertex). Thus, when the leaf $z$ is played, the weight of $z$ decreased by~$1$ and the weight of $x$ decreased by~$3$, implying that the value of $z$ is at least~$4$. However, when $z$ is played we discharge a weight of~$1$ from the value of $z$ and add a weight of~$1$ to the vertex~$x$, thereby increasing its weight to~$2$. Thus, by playing $z$ the resulting decrease in total weight is one less than the value of $z$, and is therefore at least~$4 - 1 = 3$. Further, the weight of the resulting support vertex $x$ increases from~$1$ to~$2$. We now re-color the vertex~$x$ from blue to \emph{indigo}. We note the following.

\begin{unnumbered}{Claim~\ref{t:main}.5}
If \St plays a white leaf adjacent to a purple support vertex, then the resulting decrease in total weight is at least~$3$.
\end{unnumbered}

An identical proof of Claim~\ref{t:main}.2 proves the following result.

\begin{unnumbered}{Claim~\ref{t:main}.6}
Let $R$ be the residual graph. If the game is in Phase~2 and if $C$ is an arbitrary component of $R$, then one of the following holds. \\
\indent {\rm (a)} $C \cong P_4$, with both leaves colored blue and both internal vertices colored white. \\
\indent {\rm (b)} $C \cong P_3$, with both leaves colored blue and with the central vertex colored green. \\
\indent {\rm (c)} $C \cong P_2$, with one leaf colored blue and the other colored green. \\
\indent {\rm (d)} $C \cong P_2$, with one leaf colored blue and the other colored white. \\
\indent {\rm (e)} $C \cong P_2$, with one leaf colored indigo and the other colored green.
\end{unnumbered}

A white support vertex with a white leaf neighbor in $R$ we call a \emph{targeted support vertex} in $R$. Since no two leaves are at distance~$4$ apart in $G$, every pair of targeted support vertices in $R$ are either adjacent or at distance at least~$3$ apart in $R$.

\D henceforth applies the following rules.

\begin{quote}
\textbf{Dominator's strategy:}
\begin{itemize}
\item[(R1)] \emph{Whenever \St plays a white leaf adjacent to a targeted support vertex, \D immediately responds by playing on the resulting (blue) support vertex.}
\item[(R2)] \emph{Whenever \St plays a problematic vertex, $u$, whose current value is exactly~$3$, \D immediately responds by playing the unique (targeted) support vertex associated with~$u$.}
\item[(R3)] \emph{If \D cannot play according to (R1) and (R2), he plays a targeted support vertex of maximum value.}
\item[(R4)] \emph{If \D cannot play according to (R1), (R2) and (R3), he plays a greedy move.}
\end{itemize}
\end{quote}

It remains for us to show that Dominator's strategy which applies rules (R1), (R2), (R3) and (R4) above, does indeed guarantee that on average the weight decrease resulting from each played vertex in the game is at least~$4$. We note that Dominator's strategy when playing according to (R1), (R2) and (R3) is to play a targeted support vertex or a blue support vertex with a green leaf neighbor. However, the order in which he plays such support vertices is important.
Recall that by our earlier assumptions, $G$ contains at least one leaf. The following claim will prove to be useful.

\begin{unnumbered}{Claim~\ref{t:main}.7}
While \D plays according to rule {\rm (R1)}, {\rm (R2)} and {\rm (R3)}, the following three statements hold. \\ [-2.25em]
\begin{itemize}
\item[{\rm (a)}] After each move of Dominator, every targeted support vertex has at least three white neighbors. \\ [-2em]
    \item[{\rm (b)}] After each move of Dominator, there is no green leaf adjacent to a blue vertex. \\ [-2em]
\item[{\rm (c)}] Each move that Dominator plays has value at least~$5$. \\ [-2em]
\end{itemize}
\end{unnumbered}
\proof We proceed by induction on the number, $m \ge 1$, of moves played by Dominator whenever he plays according to rule (R1), (R2) and (R3). We note that every targeted support vertex has degree at least~$3$ in the residual graph. Further, we recall that the degree sum of adjacent vertices in $G$ is at least~$4$ and the degree of a white vertex in the residual graph remains unchanged from its degree in the original graph. Since no two targeted support vertices are at distance~$2$ apart in $R$, when \D played a targeted support vertex, the white neighbors of every remaining targeted support vertex retain their color.

On Dominator's first move of the game, he plays a targeted support vertex of maximum value according to rule (R3). Such a (white) support vertex has degree at least~$3$ and all its neighbors are white, and therefore playing his first move decreases the weight by at least~$7$ and no green leaf is created. This establishes the base case when $m = 1$.  Suppose that $m \ge 2$ and that \D plays according to rule (R1), (R2) and (R3), and assume that after the first $m-1$ moves, every targeted support vertex has at least three white neighbors, there is no green leaf adjacent to a blue vertex, and each of his first $m-1$ moves has value at least~$5$. We show that after Dominator's $m$th move, the three properties (a), (b) and (c) hold.

Suppose that Staller's $(m-1)$st move plays a white leaf $x$ adjacent to a targeted support vertex $y$. Her move recolors $x$ from white to green, and recolors $y$ from white to blue. By the inductive hypothesis, before Staller played her move, the vertex $y$ had at least three white neighbors. According to rule (R1), \D immediately responds to Staller's $(m-1)$st move by playing on the resulting (blue) support vertex, $y$. Since the support vertex $y$ has at least two white neighbors after Staller played her $(m-1)$st move, his move decreases the weight by at least~$7$. Further, since the white neighbors of every remaining targeted support vertex retain their color, after Dominator's $m$th move the induction hypothesis implies that every targeted support vertex has at least three white neighbors and there is no green leaf adjacent to a blue vertex.

Suppose that Staller's $(m-1)$st move plays neither a white leaf adjacent to a targeted support vertex nor a problematic vertex. In this case, the white neighbors of every remaining targeted support vertex retain their color after Staller's move. If there remains a targeted support vertex, then, according to rule (R3),  Dominator's $m$th move plays a targeted support vertex. By induction, such a support vertex has at least three white neighbors, and therefore has value at least~$7$. Thus, as before, the desired properties (a), (b) and (c) follow by induction after Dominator's $m$th move.

Suppose that Staller's $(m-1)$st move plays a problematic vertex, $u$, whose current value is at least~$4$. Applying our discharging arguments, every targeted support vertex associated with~$u$ is recolored from white to purple. The only targeted support vertices affected by Staller's move, in the sense that it or at least one of its white neighbors changes color, are targeted support vertices associated with $u$ or adjacent to $u$. Thus, as before, the desired properties (a), (b) and (c)  follow by induction after Dominator's $m$th move.

Suppose, finally, that Staller's $(m-1)$st move plays a problematic vertex, $u$, whose current value is~$3$. In this case, either $u$ is a blue leaf with a white neighbor or $u$ is a white vertex with exactly one white neighbor. Further, the two internal vertices of a problematic path associated with~$u$ are unique. Let $uvwx$ be such a problematic path associated with~$u$, and so $v$ and $w$ are unique. In fact, $v$ is the only white neighbor of $u$. By the inductive hypothesis, immediately before Staller played her $(m-1)$st move, the targeted support vertex $w$ has at least three white neighbors. Since the vertex $v$ is the only such white neighbor of $w$ that is adjacent to $u$, after Staller plays $u$, the (white) support vertex $w$ has at least two white neighbors, including the white leaf neighbor $x$. According to rule (R2), \D immediately responds by playing as his $m$th move this unique targeted support vertex, $w$, associated with~$u$. Since $w$ has at least two white neighbors, it has value at least~$5$. As observed earlier, the white neighbors of every remaining targeted support vertex retain their color after Dominator's move. Therefore, after Dominator's $m$th move, the desired properties (a), (b) and (c) hold.~\smallqed

\medskip

By Claim~\ref{t:main}.7, while \D plays according to rule (R1), (R2) and (R3), each move he plays has value at least~$5$. By Claim~\ref{t:main}.1,  Claim~\ref{t:main}.4 and Claim~\ref{t:main}.5, each move of Staller's decreases the weight by at least~$3$. Hence, each move \D plays during this stage of the game, together with Staller's response, decreases the weight by at least~$8$. Therefore, we may assume that at some stage the game, \D cannot play according to (R1), (R2) and (R3), for otherwise  Inequality~(\ref{Eq1}) is satisfied upon completion of the game and \D can achieve his $4$-target. We note that at this stage of the game, there no longer exists a targeted support vertex. Further, there is no green leaf adjacent to a blue vertex. This implies that no green leaf adjacent to a blue vertex can be created in the remainder of the game.

According to rule (R4), \D now plays a greedy move and he continues to do so until the game is complete.  We may assume that at some stage the game enters Phase~2, for otherwise once again Inequality~(\ref{Eq1}) is satisfied upon completion of the game and \D can achieve his $4$-target. By Claim~\ref{t:main}.6 and our observation that there is no green leaf adjacent to a blue vertex when the game is in Phase~2, if $C$ is an arbitrary component of the residual graph $R$ at this stage of the game when \D cannot play according to (R1), (R2) and (R3), then $C \ncong P_2$ with one blue and one green vertex. That is, $C$ satisfies one of (a), (b), (d) or (e) in the statement of Claim~\ref{t:main}.6.

If $C \cong P_4$, then $C$  satisfies statement~(a) of Claim~\ref{t:main}.6, implying that $\w(V(C)) = 8$ and exactly two additional moves are required to totally dominate the vertices $V(C)$. If $C \cong P_3$ or if $C \cong P_2$, then $C$  satisfies statement~(b), (d) or~(e) of Claim~\ref{t:main}.6, implying that $\w(V(C)) = 4$ and exactly one move is played in $C$ to totally dominate the vertices $V(C)$. Analogously as in the proof of Claim~\ref{t:main}.3, this implies that \D can achieve his $4$-target. Thus, since $G$ has $n$ white vertices and has weight~$\w(G) = 3n$, \D can make sure that the average decrease in the weight of the residual graph resulting from each played vertex in the game is at least~$4$. Thus, in the colored-graph $G$,
$\gtd(G) \le w(G)/4 = 3n/4$.~\qed

\medskip
As an immediate consequence of the proof of Theorem~\ref{t:main} (see Claim~\ref{t:main}.3), we have the following result.

\begin{cor}
If $G$ is a colored-graph with $\delta(G) \ge 2$ and $\D$ follows a greedy strategy, then he can achieve his $4$-target.
\label{cor1}
\end{cor}

Corollary~\ref{cor1} in turn implies  Corollary~\ref{c:Tdom2}. Recall its statement.

\noindent \textbf{Corollary~\ref{c:Tdom2} {\rm (\cite{BuHeTu15})}}. \emph{The $\frac{3}{4}$-Game Total Domination Conjecture is true over the class of graphs with minimum degree at least~$2$.}

\medskip
\proof Let $G$ be a graph with $\delta(G) \ge 2$. Coloring the vertices of $G$ with the color white we produce a colored-graph in which every vertex is colored white. In particular, we note that $G$ has $n$ white vertices and has weight~$\w(G) = 3n$. By Corollary~\ref{cor1}, \D can achieve his $4$-target by following a greedy strategy. Thus, \D can make sure that the average decrease in the weight of the residual graph resulting from each played vertex in the game is at least~$4$. Thus, in the colored-graph $G$,
$\gtd(G) \le w(G)/4 = 3n/4$.~\qed

\section{Summary}

As remarked earlier, the authors in~\cite{BuHeTu15} prove a stronger result than Corollary~\ref{c:Tdom2} by showing, using game transversals in hypergraphs, that if $G$ is a graph on $n$ vertices with minimum degree at least~$2$, then $\gtd(G) \le \frac{8}{11}n$. However, our result, namely Corollary~\ref{cor1}, is surprising in that \D can complete the total domination game played in a graph with minimum degree at least~$2$ in at most~$\frac{3}{4}n$ moves by simply following a greedy strategy in the associated  colored-graph in which every vertex is initially colored white. Our main result, namely Theorem~\ref{t:main}, shows that the $\frac{3}{4}$-Game Total Domination Conjecture holds in a general graph $G$ (with no isolated vertex) if we remove the minimum degree at least~$2$ condition, but impose the weaker condition that the degree sum of adjacent vertices in $G$ is at least~$4$ and the requirement that no two leaves are at distance~$4$ apart in $G$.

\section*{Acknowledgements}

Research of both authors was supported by a grant from
the Simons Foundation (\#209654 to Douglas Rall). The first author is supported in part by the South African National Research Foundation and the University of
Johannesburg.

\end{document}